\documentclass[letterpaper, 10 pt, conference]{ieeeconf}  

\IEEEoverridecommandlockouts                              
\usepackage{authblk}
\usepackage[ruled,linesnumbered]{algorithm2e}
\DeclareMathAlphabet{\mathcal}{OMS}{cmsy}{m}{n}
\usepackage{amssymb,amsmath,bm}
\usepackage{mathtools}

\usepackage{hyperref,url}
\usepackage{comment}

\DeclarePairedDelimiterX{\norm}[1]{\lVert}{\rVert}{#1}

\usepackage{wasysym}

\usepackage[dvipsnames]{xcolor}
\usepackage{xspace}

\title{\LARGE \bf Differentiable Simulator For Dynamic \& Stochastic Optimal Gas \& Power Flows}

\vspace{-0.5cm}
\author[1]{Criston Hyett}
\author[1]{Laurent Pagnier}
\author[2]{Jean Alisse}
\author[2]{Igal Goldshtein}
\author[2]{\authorcr Lilah Saban}
\author[1]{Robert Ferrando}
\author[1]{Michael Chertkov}
\affil[1]{Program in Applied Mathematics \& Department of Mathematics, University of Arizona, Tucson, AZ, USA}
\affil[2]{Noga, The Israel Independent System Operator, Haifa, Israel}

\begin{document}

\maketitle
\pagestyle{empty}

\begin{abstract}
In many power systems, particularly those isolated from larger intercontinental grids, reliance on natural gas is crucial. This dependence becomes particularly critical during periods of volatility or scarcity in renewable energy sources, further complicated by unpredictable consumption trends. To ensure the uninterrupted operation of these isolated gas-grid systems, innovative and efficient management strategies are essential. This paper investigates the complexities of achieving synchronized, dynamic, and stochastic optimization for autonomous transmission-level gas-grid infrastructures. We introduce a novel methodology grounded in differentiable programming, which synergizes symbolic programming, a conservative numerical method for solving gas-flow partial differential equations, and automated sensitivity analysis powered by SciML/Julia. Our methodology redefines the simulation and co-optimization landscape for gas-grid systems. We demonstrate efficiency and precision of the methodology by solving a stochastic optimal gas flow problem, phrased on an open source model of Israel's gas grid model.
\end{abstract}

\section{Introduction \& Background}

The surge in renewable energy integration has heightened the variability in power demand, intensifying the fluctuations represented by the duck curve. Concurrently, the shift from coal to cleaner "bridge fuels" like natural gas places increased dependence on the gas infrastructure. This reliance extends beyond power generation to include transmission-level gas systems, which are also impacted by residential, commercial distribution, and exports. The disparate response times between gas and power networks -- seconds for power systems versus hours for gas systems -- add complexity to real-time and day-ahead coordination across these sectors. Earlier research, like that of \cite{zlotnik2016coordinated} and \cite{byeon2019unit}, integrated gas dynamics into day-ahead planning through optimization models that simplified gas network constraints. More recent efforts have developed linear approximations for pipe segments to balance computational efficiency against model fidelity, aiding their incorporation into optimization frameworks \cite{baker2023linear}. Yet, efficiently and scalably addressing the nonlinearity inherent in gas system dynamics, especially under stress and uncertainty, continues to pose a significant challenge.

The challenge we face is formally defined as solving a PDE-constrained optimization problem, which is schematically represented as:
\begin{equation}\label{eq:opt}
    \begin{aligned}
        & \min\limits_{\{u^{(s)}(t),q^{(s)}(t)\}} \sum_{s \in \mathcal{S}} \int_0^T C^{(s)}(u^{(s)}(t),q^{(s)}(t)) \, dt, \\
        & \hspace{0.2cm}\text{s.t.} \ \forall s, \ \forall t: \ \text{PDE}^{(s)}(u^{(s)}(t),q^{(s)}(t)) = 0,
    \end{aligned}
\end{equation}
where \( u^{(s)}(t) \) and \( q^{(s)}(t) \) signify the time-evolving state space and control degrees of freedom for scenarios or samples \( s \in \mathcal{S} \) respectively. The term \( C^{(s)}(u^{(s)}(t),q^{(s)}(t)) \) denotes the cumulative cost. In our chosen framework: \( q^{(s)}(t) \) embodies the gas extraction from the system, which can be redistributed across various nodes of the gas-grid where gas generators are positioned; \( u^{(s)}(t) \) represents the gas flows, gas densities, and, indirectly via the gas equation of state, pressures over the gas-grid pipes. The cost function \( C^{(s)}(u^{(s)}(t),q^{(s)}(t)) \) encapsulates the discrepancy between aggregated energy generation (directly related to gas extraction at nodes) and demand, operational costs of gas generators, and pressure constraints at the gas-grid nodes. The equation \( \text{PDE}^{(s)}(u^{(s)}(t),q^{(s)}(t))=0 \) characterizes the gas-flow equations, elucidating for each scenario \( s \) how gas flows and densities are spatially (across the gas-grid network) and temporally distributed, contingent on the profile of gas extraction and injection. A detailed explanation is provided in Section \ref{sec:meth}.

In this paper, we propose a novel approach to solving Eq.~(\ref{eq:opt}), aiming to enhance the fidelity of gas accounting in day-ahead planning of power generation in a computationally efficient manner. Our solution crafts a differentiable simulator by leveraging the principles of differentiable programming (DP)~\cite{innes2019differentiable}, combined with an efficient explicit staggered-grid method~\cite{gyrya2019explicit}, symbolic programming and the robust capabilities of the SciML sensitivity ecosystem~\cite{rackauckas2020universal}\cite{rackauckas2017differentialequations}. As we delve further, it will become evident that our approach adeptly addresses the intertwined challenges of nonlinearity, dimensionality, and stochastic modeling.

In the proposed framework, differentiable programming facilitates the calculation of gradients by seamlessly solving the gas-flow PDE across a network. This is realized by auto-generating the corresponding adjoint equations, providing flexibility in formulating the forward pass. The approach not only supports sensitivity analysis but, with a judicious selection of algorithms, proficiently manages scalability issues in parameter spaces, all while preserving the intricate nonlinear dynamics.

In this manuscript, we leverage recent efforts\cite{ma2021modelingtoolkit} in symbolic computing to decouple a model specification (i.e., PDE, boundary conditions and numerical method) and its implementation (compiled code). Historically, the user has needed to implement both the numerical method and the code representation, worrying about low-level details like parallelism, memory access, etc. These low-level details often compete with correctness of the implementation in the developer's attention. ModelingToolkit.jl is a julia package that transforms an ODE specification to compiled code - compiling and optimizing according to the target hardware (e.g. CPU/GPU). In our work, we discretize the PDE symbolically, then stitch the resulting systems of ODEs into a single system - similar to \cite{jonesMethod} but with the added detail of network topology. This system is then compiled into executable code. This approach ensures the code exactly implements the model, abstracting away from low-level implementation details. Additionally, this weakens the coupling between component pieces - allowing future users to extend this modular methodology directly, e.g. implementing compressors simply by implementing a new network component type, and writing equations describing its interaction with incident pipes.


Motivated by the everyday operational challenges characteristic of Israel's power system, as expounded in \cite{hyett2023control} and its associated references, we design and solve a dynamic, stochastic problem that integrates power and gas flows over an operational timeframe ranging from several hours to an entire day. The example provided demonstrates the following unique aspects of the system:
\begin{itemize}
    \item[(a)] Limited availability or operational restrictions of gas compressors;
    \item[(b)] Notable fluctuations in renewable resources and power loads, with curtailment being inadmissible under the normal operational paradigms assumed in this research;
    \item[(c)] An intentionally over-engineered power system, ensuring power lines remain within thermal boundaries during standard operations.
\end{itemize}
Note critically that while the example presented in the following contains gas network specializations (e.g., (a)), because of the generality of symbolic programming and automatic differentiation, the proposed methodology is not restricted to these simplifying assumptions on the gas network.

The remainder of the manuscript is structured as follows: In Section \ref{sec:meth}, we elucidate our gas modeling methodology, elaborate our fundamental optimization problem, and delineate our strategy for its resolution. Experimental results for a representative regional gas network are presented in Section \ref{sec:res}. Finally, the manuscript culminates with conclusions and suggested future directions in Section \ref{sec:con}.

\section{Methodology}\label{sec:meth}
\subsection{Solving PDE Constrained Optimization}

\begin{figure}
    \centering
\includegraphics[width=\linewidth]{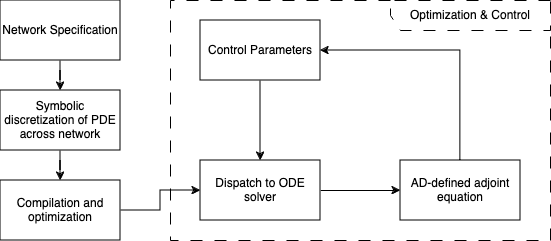}
    \caption{Schematic of methodology}
    \label{fig:methodology}
\end{figure}

In this Section, we elucidate our strategy to address Eq.~(\ref{eq:opt}). Essentially, two predominant methodologies emerge for tackling the PDE-constrained optimization challenge: 
\begin{enumerate}
    \item \textbf{Constraint Matrix Encoding:} This method integrates the PDE into a constraint matrix that grows as discretization becomes finer. A notable merit of this approach is its flexibility in harnessing advanced optimization techniques and the ability to consider all timepoints simultaneously. However, the methodology grapples with potential pitfalls such as the emergence of unphysical solutions, non-adherence to constraints during intermediary timeframes, and the curse of dimensionality, manifesting as an exponential surge in complexity with the growth of the problem.

    \item \textbf{Differentiable Programming:} This strategy leverages the adjoint method (standard material included for completeness in Appendix~\ref{app:adj}) to calculate gradients with respect to control parameters \textit{through the PDE solver}. This method ensures the PDE solution $u$ remains physically valid throughout the optimization process. Further, $u$ converges to a well-defined solution as the grid undergoes refinement ($\Delta x,\Delta t \to 0$). However, challenges arise in generalizing the calculation of gradients through the PDE solver, dimensionality of the discretized PDE growing as the grid refines, and induced temporal computational complexity due to the Courant-Friedrichs-Lewy (CFL) condition for hyperbolic PDEs\cite{brio2010numerical}.

\end{enumerate}

In the present study, we adopt the second approach and implement it as outlined in Fig.~(\ref{fig:methodology}). We confront the aforementioned challenges by 
\begin{enumerate}
    \item \textbf{Using symbolic programming to discretize the PDE, and pre-compute coupling conditions at nodes.} This eliminates unnecessary memory accesses during simulation, ensures compatibility across network topologies, and creates an auto-differentiation friendly system of ODEs. Further, efficient memory access and the composability of julia enable easy parallelization across threads, processors, or accelerators\cite{rackauckas2017differentialequations}.
    
    \item \textbf{Using a conservative, method-of-lines discretization to allow for temporal integration using high-order, strong-stability preserving numerical integrators} to mitigate temporal computational complexity induced by the hyperbolic structure of the PDE.
    \item \textbf{Leveraging modern, source-to-source Automatic Differentiation (AD) tools to automatically define and solve the corresponding adjoint equation}; allowing for freedom in expanding the network component library, and favorable scaling when computing gradients of high-dimensional parameterizations.\cite{innes2019differentiable}
\end{enumerate}


Collectively, our proposed methodology bridges the gap between low-fidelity, optimization-centric 'constraint matrix' methods suited for long-term planning, and the demand for a physics-based tool tailored for medium-term to real-time planning and analysis, essential for operational coordination between power and gas utilities.

\subsection{Gas-Flow Equations}

We begin by discussing the dynamics of a single pipe. The governing partial differential equations (PDEs) for the Gas Flow (GF), describing the dynamics of density $\rho(x,t)$ and mass-flux $\phi(x,t)$ along the pipe coordinate $x$ with respect to time $t$, are provided as follows \cite{osiadacz1984simulation},\cite{steinbach_pde_2007},\cite{chaudhry2014applied}:
\begin{align}
    & \partial_t \rho + \partial_x \phi = 0, \label{eq:gf_mass} \\
    & \partial_t \phi + \partial_x p = -\frac{\lambda}{2D} \frac{\phi |\phi|}{\rho}, \label{eq:gf_momentum}
\end{align}
where $\lambda$ is the Darcy-Weisbach friction factor.

These equations are valid under the assumption that the gas velocity is much smaller than the speed of sound in the gas ($\phi/\rho \ll a$). This is a reasonable approximation for the typical flows we consider.

To provide a complete description, it is necessary to relate the pressure $p(x,t)$ and density $\rho(x,t)$ using an equation of state:
\begin{equation}\label{eq:eos}
    p = Z(\rho,T)RT\rho,
\end{equation}
where $Z(\rho, T)$ denotes the compressibility factor. For clarity, we adopt the ideal gas law to model the equation of state, where $Z(\rho,T)RT$ is replaced by a constant, $a^2$, with $a$ representing the speed of sound in the gas. Notably, there are more accurate models available (e.g., CNGA \cite{menon}), and the methodology we present here is agnostic to the specific choice of model.

The system of Eqs.~(\ref{eq:gf_mass},\ref{eq:gf_momentum},\ref{eq:eos}) is also supplemented by the boundary conditions, for example given profile of injection/consumption, at both ends of the pipe of length $l$, 
\begin{equation} \label{eq:flow_bcs}
    q(0,t) \text{ and } q(L,t).
\end{equation}

To extend the described equations from a single pipe to a network, the boundary condtions (\ref{eq:flow_bcs}) need to appropriately couple pipe boundaries together, depending on the network topology. We will discuss it below.


\subsection{Explicit Numerical Method for the Forward Path}

To solve Eqs.~(\ref{eq:gf_mass},\ref{eq:gf_momentum},\ref{eq:eos}) in the case of a single pipe and their network generalizations, we use a conservative, explicit, staggered-in-space-grid method motivated by the staggered-grid (in space and time) method of Gyrya \& Zlotnik \cite{gyrya2019explicit} but differentiated by rejecting the staggering in time. As discussed later, much of the computational complexity in the simulation of this PDE is from long time simulation while restricted to small timesteps. Rejecting the staggering in time allows for deploying the discretized system to higher-order adaptive methods.
This yields a method-of-lines style discretization that can be solved with a strong stability preserving (SSP) ODE integrator.
\begin{figure}
    \centering
    \includegraphics[width=\linewidth]{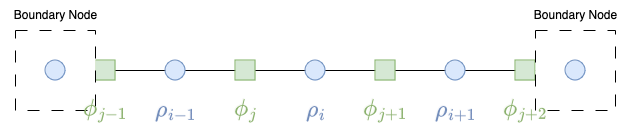}
    \caption{Schematic of spatially staggered grid with variable locations}
    \label{fig:staggered_grid}
\end{figure}
As illustrated in Fig.~(\ref{fig:staggered_grid}), we let the staggered grids for $\rho$ and $\phi$ be denoted by $x_i, x_j$ respectively, such that $x_j = x_i+\Delta x/2$. Futher, denote $\rho^n_i := \rho(x_i, t_n)$ and $\phi^n_j := \phi(x_j, t_n)$. Then we take centered differences on the staggered grid of Eqs.~(\ref{eq:gf_mass},\ref{eq:gf_momentum}) to obtain
\begin{align} \label{eq:num_rho_update}
    \frac{d\rho^n_i}{dt} &= - \frac{1}{\Delta x}\left(\phi_{j}^n - \phi_{j-1}^n\right),\\ \label{eq:num_phi_update}
    \frac{d\phi_j^{n}}{dt} &= - \left(\frac{p^{n}_{i+1} - p^{n}_i}{\Delta x}+\frac{\lambda}{2D} \frac{\phi^{n}_j|\phi^{n}_j|}{\rho_j^{n}}\right).
\end{align}
Here $\rho_j^n \approx \frac{1}{2}(\rho^n_{i-1} + \rho^n_{i})$.
\begin{figure}
    \centering
    \includegraphics[width=0.7\linewidth]{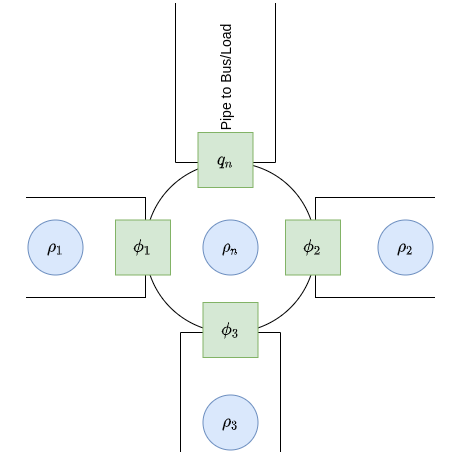}
    \caption{An example of a 3-pipe junction discretization, with a fictitious pipe setting the nodal flow boundary condition.}
    \label{fig:num_boundary}
\end{figure}
As we are interested in integration in day-ahead planning of energy generation, we control Dirichlet boundary conditions on nodal mass flows $\{q_i\}_{i \in \text{nodes}}$ (directly relating to generated power through heat-rate curves). These boundaries are resolved according to the numerical method using a boundary discretization shown in Fig.~(\ref{fig:num_boundary}). The density updates for these junctions are evaluated using conservation of mass at the boundary node $\ell$
\begin{equation} \label{eq:num_bc}
\begin{split}
\frac{d\rho_\ell^n}{dt}\sum_{k \in \partial \ell} &\left( \frac{\Delta x}{2} \right)S_{kl} =\\ &q_\ell + \sum_{k \in \partial \ell} \text{sgn}_{kl}S_{kl}\phi_{kl}^n(x=x_\ell),
\end{split}
\end{equation}
where $S_{kl}$ is the cross-section area of pipe from node $k$ to node $l$, and $\text{sgn}_{kl}$ keeps track of the directionality of the mass flux. $\rho_{kl}, \phi_{kl}$ denote the $l$-side boundary values of density and mass flux for the pipe from node $k$ to node $l$. $\phi_{kl}^n(x=x_\ell)$ is approximated by a second-order, one-sided stencil.
After solving for the density at the node, the flux update at the ends of the pipes can proceed using the momentum equation (\ref{eq:num_phi_update}).

\begin{figure}
    \centering
    \includegraphics[width=0.8\linewidth]{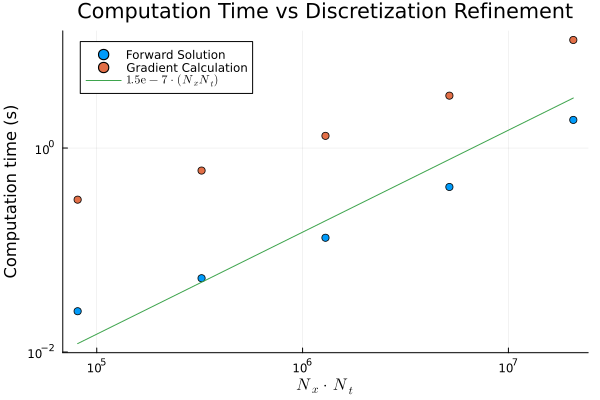}
    \includegraphics[width=0.8\linewidth]{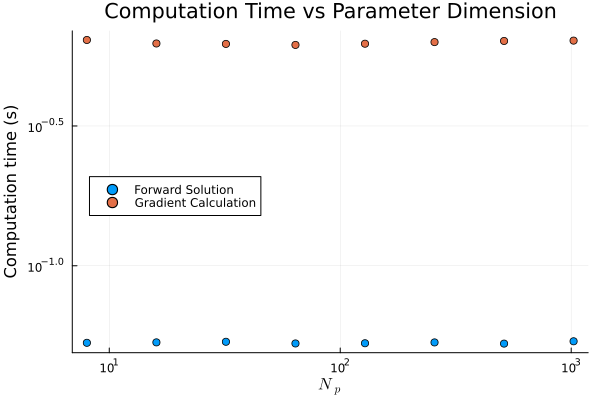}
    \caption{Computational complexity scaling for the forward and adjoint calculations, (top) as a function of dimensionality of discretized PDE, and (bottom) as a function of parameter dimension.}
    \label{fig:computation_complexity}
\end{figure}


\subsection{Optimization Formulation: Cost Function}

In our pursuit to devise a scalable framework that aptly accommodates optimization challenges akin to the archetype presented in Eq.~(\ref{eq:opt}), we pivot our attention to a paradigmatic problem: the minimization of an integrated objective spanning time and evaluated under the cloak of uncertainty. This uncertainty, reflected through diverse scenarios \(s \in \mathcal{S}\), pertains to the gas injection consumption \(q^{(s)}(t):=\{q_i^{(s)}(t)\}_{\forall i \in \text{nodes}}\), influenced possibly by variable renewable generation. The time interval \(t\in[0,T]\) typically encapsulates a pre-established performance window.

Our control parameters, symbolized by nodal flows \(q^{(s)}(t)\), permit adjustments within our forthcoming dynamic and stochastic optimization context. 
Letting \(u^{(s)}(t) := \{\rho^{(s)}(x,t), \phi^{(s)}(x,t)\}_{\forall x \in \text{nodes \& pipes}}\) be the solution of the PDE defined by Eqs.~(\ref{eq:num_rho_update}-\ref{eq:num_bc}), with nodal flows \(q^{(s)}(t)\).
\begin{figure}
    \centering
    \includegraphics[width=\linewidth]{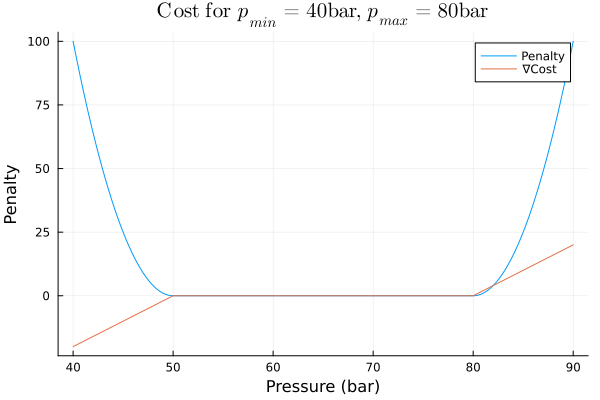}
    \caption{Quasi-quadratic penalty for violating pressure constraints, with $40-80$ bar shown here, but configurable on a per-node basis.}
    \label{fig:opt_form_pressure_penalty}
\end{figure}

Our primary optimization task is delineated as minimization of
\begin{equation} \label{eq:big_pic_opt}
O(u,q) =\sum_{s \in \mathcal{S}} \int_0^T C(u^{(s)}(t),q^{(s)}(t)) \, dt,
\end{equation}
where the specific per-time and per-scenario cost is expanded as:
\begin{equation} \label{eq:method_objective}
\begin{aligned}
    C(u^{(s)}&(t), q^{(s)}(t)) = \alpha \left(D(t) - \sum_{i\in \text{nodes}} G_i(q^{(s)}_i(t))\right)^2\\
    & + \beta \sum_{i \in \text{nodes}} E_i(q^{(s)}_i(t)) + \gamma \sum_{x\in \text{nodes}} V(p^{(s)}(x,t)),
\end{aligned}
\end{equation}
constrained by the gas-flow PDEs and associated boundary conditions over gas-grid network detailed earlier. 

The first term in Eq.~(\ref{eq:method_objective}) aims to minimize the cumulative mismatch between energy demand $D^{(s)}(t)$ and the sum of generation at each node $i$ and at each moment of time $t$, $G_i(q^{(s)}_i(t))$, with $q_i^{(s)}(t)$ representing the nodal flows, which is our control variable (one we are optimizing over). $G_i(q^{(s)}_i(t))$ is an efficiency function, mapping mass flow (in $kg/s$) to power production (in $MW$). 
Here the assumption is that any residual mismatch, if not optimal, can be adjusted by either shedding demand or introducing a generation reserve, at a certain cost.
The second term in Eq.~(\ref{eq:method_objective}), \(E_i(q_i(t))\), stands for the cost of operating power generator run on gas and located at the node \(i\) at the gas withdrawal rate \(q_i(t)\). The third term in Eq.~(\ref{eq:method_objective}), \(\sum_{x\in \text{nodes}} V(p^{(s)}(x,t))\), is chosen to be a quasi-quadratic cost (regularized by the $\text{relu}$ function) to penalize pressure constraint violations across the network (refer to Fig.~\ref{fig:opt_form_pressure_penalty}): with \(p_{\text{min}}(x)\) and \(p_{\text{max}}(x)\) denoting pre-set pressure boundaries at system nodes. The influence of the multi-objective cost \(C\)'s components can be modulated using the hyperparameters \(\alpha\), \(\beta\), and \(\gamma\).

\section{Results}\label{sec:res}

\begin{figure}
    \centering
    \includegraphics[height=2.5in]{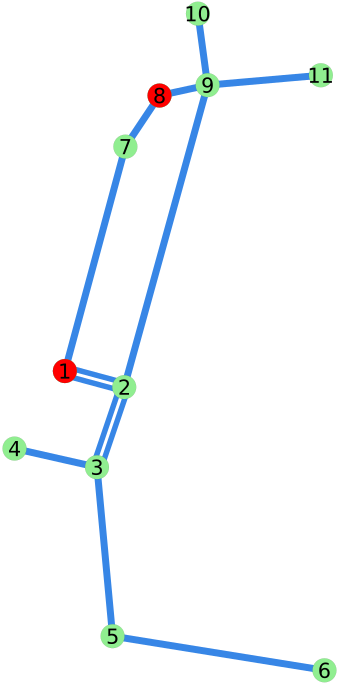}
    \caption{11-node network diagram, with supply nodes (injection into the network) in red, demand nodes (consumption from the network) in green.}
    \label{fig:11_pipe_opt}
\end{figure}


To exemplify the methodology, we solve an optimal gas flow problem phrased on a previously studied reduced model of Israel's gas grid \cite{hyett2023control}. The reduced model, shown in Fig.~(\ref{fig:11_pipe_opt}), has 11 nodes, with a total pipe length of approximately $550$km. We minimize the objective in Eq.~(\ref{eq:big_pic_opt}) over a time horizon of 10hrs, encompassing a morning ramp in energy use. The demand curves $D(t)$ are aggregated from publicly available data; the gas cost $G(q)$ is taken as a constant; the efficiency curves $E_i(q(t))$ take one of three constant values representing efficient, nominal, and inefficient turbines; and the pressure limits are set as $p_{min} = 60\text{bar}$, $p_{max} = 80\text{bar}$. We use a box-constrained Limited-memory Broyden–Fletcher–Goldfarb–Shanno (LBFGS) optimizer, the constraints enforcing max and min injection/consumption at each node. Each node has hourly flow-rate control parameters, so the dimension of the optimization space is $110 = 11\text{nodes} \ * 10\text{hrs}$.

\begin{figure}
    \centering
    \includegraphics[width=0.8\linewidth]{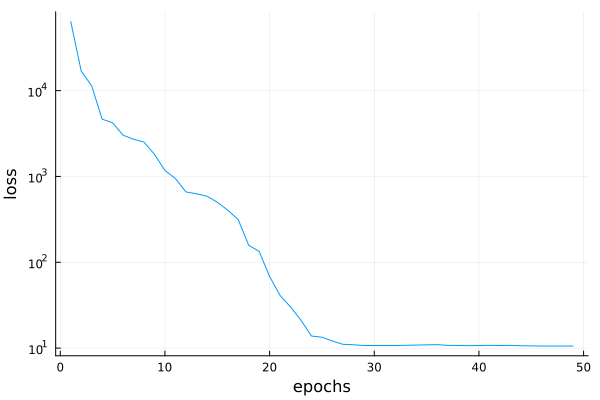}
    \vspace{0.1in}
    \includegraphics[width=0.8\linewidth]{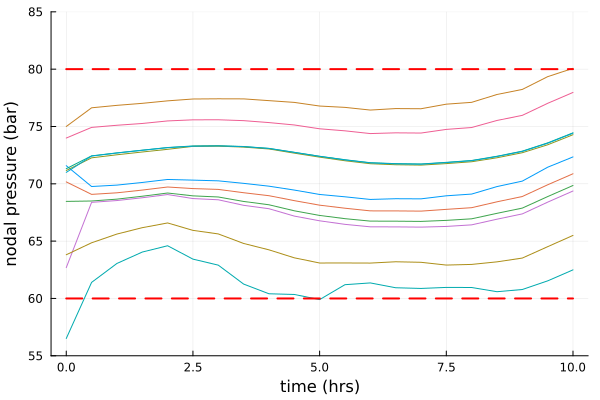}
    \vspace{0.1in}
    \includegraphics[width=0.8\linewidth]{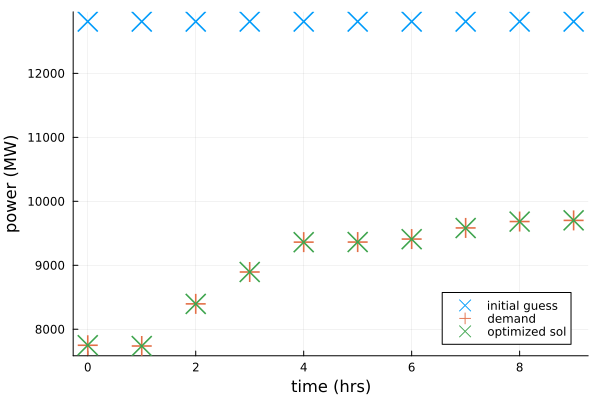}

    \caption{Optimal gas flow results from an 11-node network, spanning $\approx 550$km. The OGF is solved across a time horizon of 10 hours, using representative data of a morning ramp in the energy demand. (Top) Shows the quick convergence of LBGFS, despite the initial guess yielding a large penalty. (Middle) shows the evolution of pressure at each node using the optimized injections/withdrawals. Note that despite the dynamic initial conditions being outside the pressure window, the optimization quickly rectifies and holds all pressures in the acceptable range marked by the dashed lines. (Bottom) shows that we meet demand during the morning ramp, without waste.}
    \label{fig:ogf}
\end{figure}


The optimization results, considering a deterministic gas consumption profile, are depicted in Fig.~\ref{fig:ogf}. We observe an exponential decrease in the loss function, with the algorithm converging to a stable minimum within $30$ iterations. System pressures remain within specified limits, except for the uncontrolled initial data, which does not contribute to the loss calculation. Notably, at node 6 -- the system's lowest pressure point -- there is a proactive pressure increase to accommodate the expected rise in consumption during the morning peak. This optimization strategy, navigating the system's nonlinearities, proves crucial for operators in making informed real-time decisions. Impressively, the optimizer achieves exact demand fulfillment, even starting from a sub-optimal initial condition.

\begin{figure*}
    \centering
    \includegraphics[width=0.7\linewidth]{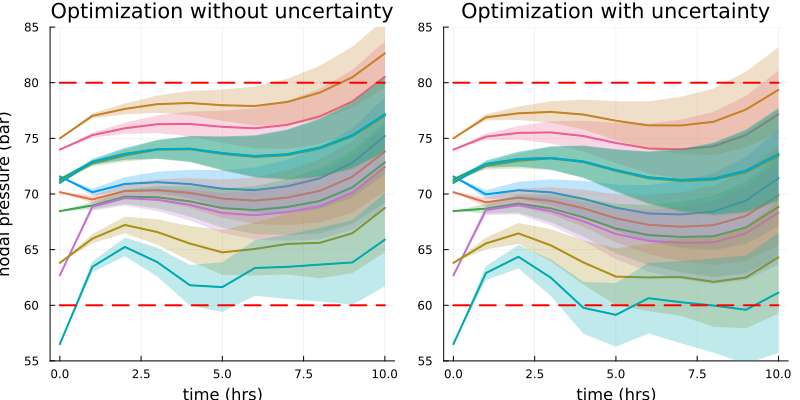}
    \caption{Results of optimization under uncertainty. We initialize the optimization using the solution in Fig.~(\ref{fig:ogf}), and then re-perform the optimization, taking the loss as the expectation over samples of noisy consumption as in Eq.~(\ref{eq:big_pic_opt}). Plotted are the mean (solid line) plus/minus a standard deviation (opaque region) for selected nodes.
    This procedure illustrates the ability to tune for robustness against stochastic fluctuations in the power network.}
    \label{fig:ogf_uncertainty}
\end{figure*}


Subsequently, we applied the software to perform optimization under uncertain consumption patterns, assuming a normal distribution with a standard deviation equal to 5\% of the current consumption level. This scenario aims to simulate the uniform response of all generators to the variability inherent in renewable energy sources. The findings, illustrated in Fig.~(\ref{fig:ogf_uncertainty}), indicate the network's ability to achieve a reduced minimum. This improvement is attributed to effectively managing the pressures at nodes $8$ and $7$ (the nodes with the highest pressure) to remain below 80bar towards the latter part of the simulation, thereby mitigating less frequent low-pressure breaches at node $6$, the node with the lowest pressure.

All code required to run and analyze these simulations as well as smaller test networks is available at \url{https://github.com/cmhyett/DiffGasNetworks}.


\section{Conclusion \& Path Forward}\label{sec:con}

The chief technical contribution of this manuscript is the integration of symbolic programming, automatic differentiation, and gas-flow PDE numerics to develop a more accurate, physics-based approach for addressing optimization and control issues in gas networks. This was demonstrated through the solution of a stochastic optimal gas flow problem. Our development efforts concentrated on:
\begin{enumerate}
    \item \textbf{Efficiency:} We prioritized the forward solution's efficiency and the gradient computation's scalability using the adjoint method. Achievements in efficiency resulted from combining symbolic programming, high-performance ODE integrators, and advanced AD tools.

    \item \textbf{Consistency:} For PDE-constrained optimization, especially in short-term or real-time planning, it is crucial to maintain precise physical solutions irrespective of grid refinement choices. The differentiable programming framework ensures consistent convergence within the physical domain and offers error assurances.

    \item \textbf{Flexibility:} The methodology's design allows adaptable network and component configurations, supporting a range of applications from uncertainty quantification to inverse problems and data assimilation. By preserving symbolic representations until execution and employing Automatic Differentiation for derivative calculations, our approach facilitates selecting the most appropriate gradient computation method from available library of options, ensuring both application breadth and solution specificity.
\end{enumerate}


We demonstrated our ability to utilize these characteristics to effectively solve an Optimal Gas Flow (OGF) problem under uncertainty, preserving essential system properties and the complete nonlinear dynamics of a representative regional gas network.


Future endeavors will focus on integrating this methodology into the broader scope of gas network optimization and control. Although our current model omits gas compressors due to specific characteristics of the Israel's system, the extension of this method to include compressors and valves in the network library is straightforward. This adaptability is facilitated by the inherent generality of our method's design.


Moreover, while our Optimal Gas Flow (OGF) model under uncertainty effectively managed the expected cost, the realm of stochastic optimal control offers the capability to target more specific objectives, such as managing the higher moments of nodal pressures. Proactively addressing and planning for these uncommon occurrences within gas-grid system coordination remains a dynamic and critical field of research.


Finally, there have been significant theoretical advancements in multi-fidelity methods for outer-loop applications\cite{peherstorfer2018survey}. By integrating the methodology presented in this manuscript with constraint-matrix and machine learning techniques, we could develop a comprehensive multi-fidelity approach. This integration promises to be both efficient and versatile, offering high-fidelity solutions to the underlying Partial Differential Equations (PDEs) across networks. Such an approach has the potential for broader applicability, extending beyond gas networks to encompass a wider spectrum of complex systems.

\section*{Appendices}\label{sec:app}
\subsection{Adjoint Method} \label{app:adj}

In order to utilize gradient descent algorithms to optimize Eq~(\ref{eq:big_pic_opt}), it is necessary to compute $\nabla_q O(u, q)$, where the objective function $O(u, q)$ is defined as
\begin{equation*}
    O(u, q) = \int_0^T C(u, q, t) \, dt,
\end{equation*}
with $C$ representing the cost. Given that $g(u, \dot{u}, q, t) = 0$ and $h(u(0), q) = 0$ define the differential equation and the initial condition, respectively, we can reformulate this optimization problem using the Lagrangian formulation:
\begin{equation*}
\mathcal{L} := \int_0^T \left[ C(u,q,t) + \lambda^T g(u,\dot{u},q,t) \right]dt + \mu^T h(u(0),q),
\end{equation*}
where $\lambda(t)$ and $\mu$ being the Lagrangian multipliers. 

To compute the gradient of the Lagrangian, $\nabla_q \mathcal{L}$, we integrate by parts to express $\nabla_q \dot{u}$ in terms of $\nabla_q u$, leading to the following expression after substitutions and rearrangements:
\begin{multline}
    \nabla_q \mathcal{L} = 
    \int_0^T \Big[ \left(\partial_u C \mkern-5mu + \mkern-5mu \lambda^T \mkern-7mu \left( \partial_u g \mkern-5mu - \mkern-5mu d_t \partial_{\dot u}g \right) \mkern-5mu - \mkern-5mu \dot \lambda^T \partial_{\dot u} g \right)\nabla_q u +\\ \partial_q C + \lambda^T \partial_q g \Big]dt 
    + \left.\lambda^T \partial_{\dot u} g \nabla_q u \right|_T \\ + \left.\left( \mu^T \partial_u h - \lambda^T \partial_{\dot u} g \right)\right|_0 \nabla_q u(0) + \mu^T \nabla_q h.
\end{multline}

Resolving zero gradient conditions at $t=0$ and $t=T$ and solving for $\lambda$ and $\mu$ appropriately, we eliminate the need to directly calculate $\nabla_q u$.  This results in 
\begin{multline}
    \nabla_q \mathcal{L} = \int_0^T \Big[ \left(\nabla_q C + \lambda^T\left( \partial_u g - d_t \partial_{\dot u} g\right) + \dot \lambda^T \partial_{\dot u} g\right)\nabla_q u\\
    + \lambda^T \partial_q g \Big] dt
    + \left[ \lambda^T \partial_{\dot u}g \right|_0 \left( \partial_{u(0)} h \right)^{-1} \nabla_q h.
\end{multline}

We then solve for $\lambda$ backward in time from the differential equation involving $\lambda$, $\partial_q C$, and the derivatives of $g$ with respect to $u$ and $\dot{u}$, ensuring $\lambda(T) = 0$:
\begin{equation} \label{eq:adj_diffeq}
\begin{split}
    \partial_q C + \lambda^T \left(\partial_u g - d_t \partial_{\dot u} g\right) + \dot \lambda^T \partial_{\dot u} g = 0\\
    \text{with } \lambda(T) = 0.
\end{split}
\end{equation}

This process yields the final gradient expression for the optimization step
\begin{align} \nonumber 
    \nabla_q \mathcal{L} & = \nabla_q O =
    \int_0^T \big( \nabla_q C + \lambda^T \partial_q g \big) \, dt \\ & + \left. \dot \lambda^T \partial_{\dot u} g \right|_0 \left(\partial_{u(0)}h\right)^{-1} \nabla_q h, \label{eq:adj_int}
\end{align}
which uses Eq.~(\ref{eq:adj_diffeq}) for efficient evaluation. 

Notice, that the functional forms such as $\nabla_q C(u, q, t)$ and $\partial_u g(u, \dot{u}, q, t)$, dependent on the solved state $u$, are determined and evaluated using source-to-source Automatic Differentiation (AD).

\subsection{Differentiable Programming}

Source-to-source differentiation, particularly from Zygote.jl \cite{innes2019dont}, is a transformational capability that allows reverse-mode automatic differentiation (AD) through programming language constructs -- enabling optimized adjoint function evaluation without the need to write the derivatives by hand. This freedom ensures correctness, and allows for generality in construction of the forward pass \cite{rackauckas2022scimlbook}.

In order to compute the integral Eq.~(\ref{eq:adj_int}), the adjoint ODE Eq.~(\ref{eq:adj_diffeq}) is solved for $\lambda(t)$, and the term $\partial_q g$ is found via source-to-source reverse-mode AD. This method to compute the adjoint has computational cost that scales linearly with the forward pass, and with the number of parameters \cite{ma2021comparison}. 
Notably, the sensitivity backend used evaluates the most efficient approach, and dispatches the appropriate sensitivity method (e.g., forward or adjoint sensitivities). We emphasize the adjoint method here, because of its preferable scaling with respect to high-dimensional parameterizations.
\begin{figure}
    \centering
    \ifdefined\ifCLASSOPTIONtwocolumn
        \includegraphics[width=0.8\linewidth]{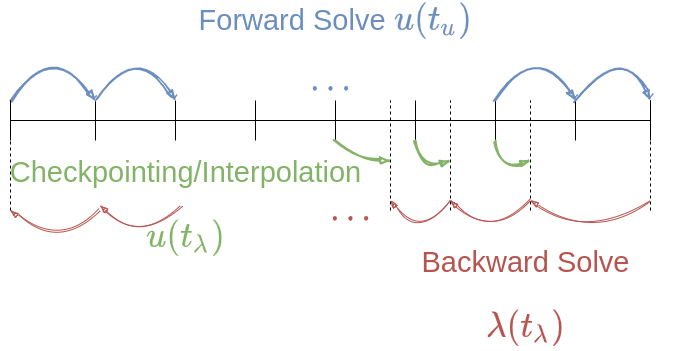}
    \else
        \includegraphics[width=0.8\linewidth]{adjoint_method.png}
    \fi
    \caption{Schematic of adjoint method implementation, where the adjoint ODE utilizes checkpointing to accelerate queries of the state $u$ at required times.\cite{ma2021comparison}}
    \label{fig:dp_adjoint}
\end{figure}


\newpage



%
\bibliographystyle{IEEEtran}
\bibliography{refs}

\end{document}